\documentclass{article}

\usepackage{fullpage}
\usepackage{graphicx}

\newtheorem{prop}{Proposition}

\usepackage{amsfonts,amsmath}

\def\Rset{\mathbb{R}}
\def\sgn{\mbox{sgn}}
\def\sat{\mbox{sat}}

\begin{document}

\title{Extremum seeking via continuation techniques for optimizing biogas
production in the chemostat}
\author{{\sc A. Rapaport$^{1}$, J. Sieber$^{2}$, S. Rodrigues$^{3}$ and M. Desroches$^{4}$}\\[2mm]
$^{1}$ UMR INRA/SupAgro MISTEA and EPI INRA/INRIA MODEMIC,
  Montpellier, France.\\
Email:~{\tt rapaport@supagro.inra.fr}\\
 $^{2}$ CEMPS, Univ. of Exeter, U.K. Email:~{\tt
   J.Sieber@exeter.ac.uk}\\
 $^{3}$ CN-CR, Univ. of Plymouth, U.K. Email:~{\tt
   serafim.rodrigues@gmail.com}\\
 $^{4}$ EPI SISYPHE, INRIA Rocquencourt, France. Email:~{\tt
   Mathieu.Desroches@inria.fr}
}

\maketitle

\begin{abstract}
We consider the chemostat model with the substrate concentration as
the single measurement. We propose a control strategy that drives the
system at a steady state maximizing the gas production without the
knowledge of the specific growth rate.
  Our approach separates the extremum seeking problem from the
  feedback control problem such that each of the two subproblems can
  be solved with relatively simple algorithms. We are then free to
  choose any numerical optimization algorithm. We give a demonstration
  for two choices: one is based on slow-fast dynamics and numerical
  continuation, the other is a combination of golden-section and
  Newton iteration.  The method copes with non-monotonic growth
functions.
\end{abstract}

\section{Introduction}

The control design of chemostat models have been extensively studied in the
literature, with the objective to provide efficient and reliable
control strategies for industrial applications, such as
biotechnological or pharmaceutical processes \cite{BD90,BV95,DP97}.
A common task is to drive a continuous stirred tank bioreactor to a
set-point that optimizes a objective function, for instance the
gas production rate \cite{DGK72,SCNBV04}. In this work, we consider the chemostat
model with a single strain \cite{SW95}
\begin{equation}
\label{chemostat}
\begin{array}{lll}
\dot s & = & -\mu(s)b+u(s_\mathrm{in}-s)\\
\dot b & = & \mu(s)b-ub
\end{array}
\end{equation}
where the state variables $s$ and $b$ denote the substrate and biomass
concentrations, respectively (in these equations the yield coefficient
has be chosen equal to one without any loss of generality). The
concentration of substrate in the feed is denoted by $s_\mathrm{in}$, and the
dilution rate $u\geq 0$ is the manipulated variable. The production
rate to be maximized is given by
\[
r=\mu(s)b \ .
\] 
In such bioprocesses, if
often happens that the growth kinetics $\mu(\cdot)$ is unknown (or
poorly known) and possibly evolves slowly with time or changes of environment (temperature,
pH...). Consequently,
the robustness of the control strategies with respect to uncertainties
on $\mu(\cdot)$
is a crucial issue for real applications. Many works have been done for the
on-line estimation of the growth function \cite{H82,HR82,DB86,VCV97,VCV98,PFFD00,BG03} and the robust stabilization of such processes about a given reference point \cite{DB84,JBV99,SC99,RH02,SAL12},
but there are comparatively much less works that concern the
maximization of an objective under model uncertainties \cite{DPG10}. Such issues
are typically addressed by the design techniques of extremum
seeking controls \cite{W64,G74,KW00,K00,AK03,AK04,TNM08,TNMA09}. Roughly speaking, these
techniques consist in adding an excitation signal to the input 
$u=\bar u +a\sin(\omega t)$, of small amplitude $a$ and high frequency
$\omega$, and capturing an estimation of
the gradient of the objective function by filtering. Three time scales
(excitation signal, process, filtering) then operate in the overall closed loop system.
Such a technique has been
developed first in \cite{WKB99} for unknown functions $\mu(\cdot)$ of Monod or
Haldane type. Later on, an adaptive extremum-seeking scheme has
been proposed to improve the performances of the transient response,
for the Monod's case in \cite{ZGD03} and the Haldane's one in
\cite{MGZD04}. 
In these approaches, the production rate $r$ is assumed to be measured
on-line, and explicit expressions of the unknown function
$\mu(\cdot)$ are required. An alternative extremum-seeking scheme 
combined with a neural network has been proposed in \cite{GDP04} to
get out of these requirements.
Nevertheless, all the mentioned approaches require the on-line measurement of the
whole state $(s,b)$, excepted in \cite{MGD04} where $b$ only is
measured but the result is dedicated to the Monod's kinetics.
In the present paper, we consider that the single on-line measurement
\begin{equation}
y=s
\label{eq:output}
\end{equation}
is available. One may consider coupling the former extremum-seeking
controllers with state observers
\cite{DPY92,BD90,GHO92,D03,RM04,RD04,MD07,KJ13}, but there is here an intrinsic
difficulty due to the fact that ``good'' observers, i.e. whose speed of
convergence can be arbitrary tuned, require the
exact expression of the dynamics, and that ``fast'' observers are
often sensitive to measurement noise and input disturbances.

In the present work, we present a rather different approach that
does not require any knowledge on the growth function neither the
measurement of the objective function, and that is based on a continuation
method \cite{AG03}. We do not consider any excitation signal. Instead, we design a
dynamical extension of the chemostat model (\ref{chemostat}) in a
slow-fast or ``singularly perturbed'' form \cite{O91}, such
that the attracting critical manifold is exactly the graph of the unknown
function $s\mapsto \mu(s)$. 

We present here a rather different and new approach that does
not require any knowledge of the growth function. Neither do we
require the ability to measure the objective function at all
times. Instead of feeding an excitation signal into the plant (the
usual approach in extremeum seeking) we apply a classical feedback
control stabilizing the plant to guarantee a uniform decay of
transients. This stabilization step is then combined with classical
numerical methods that treat the equilibrium output of the
feedback-stabilized plant as an unknown function of the input. We
demonstrate this approach with two different numerical methods. 
We apply a continuation method \cite{AG03} tracing out the graph of equilibrium
outputs (depending on the inputs) until we reach the maximum of the
objective. We then apply a classical optimization method to the graph,
in our case a golden section iteration and Newton iterations.

Our approach splits the dynamic extrmemum seeking
problem into two subproblems: a feedback stabilization problem and a
numerical optimization problem. This has several advantages:
\begin{itemize}
\item First, the proof of convergence of the overall scheme is split into
simpler subproblems: stability of the feedback for fixed inputs and
convergence of the numerical algorithm. In the case of the chemostat
with a growth rate of Haldane type the first part (stability) is
globally true \cite{SRRD13}. The second part (convergence of
numerics) follows in general at least locally from the general
theorems developed for numerical algorithms. For the specific case of
the chemostat with growth rate of Haldane type one can again state
global convergence for the continuation and the golden-section
iteration because of the uniqueness of the local maximum.
\item Second, the number of different time scales is reduced to two: the
decay time of the feedback controlled system and the iteration process
of the numerical algorithm. In theory the decay time can be made
arbitrarily short by adjusting the feedback gains. In practice this is
limited by intrinsic noise, measurement noise and sampling time. The
numerical convergence speed can be superlinear (quadratic) in the
ideal case for Newton iterations.
\item Third, one only needs to evaluate the objective at equilibrium. One
can see the importance of this aspect in the chemostat example: the
rate $r=\mu(s)b$ depends on the unknown $\mu$ and the internal state
$b$, which is not measured. However, in equilibrium $\mu$ equals the
input $u$ and $b$ equals $s_\mathrm{in}-s$ such that we can apply the
continuation or optimization to $D(s_\mathrm{in}-s)$ instead.
\end{itemize}
Our extremum seeking scheme is an extension of a
former work \cite{SRRD12,SRRD13} that gives such a method for the reconstruction of the
graph of $\mu(\cdot)$ with the single measurement of the substrate, without any a priori
knowledge on $\mu$ (excepted to be a smooth function). In particular,
this technique can cope with non-monotonic growth functions
and allows to identify unstable states of the open-loop system. The
paper is organized as follows. In Section \ref{section_general} we
first present our general methodology, and then show how to apply it on the
chemostat model in Section \ref{section_chemostat}. Section
\ref{section_demo} is devoted to numerical demonstrations.

\section{Assumptions in the general framework}
\label{section_general}
We consider a single-input/single-output dynamical system
\begin{equation}
\begin{array}{lll}
\dot x & = & f(x,u)\\
y & = & h(x)
\end{array}
\end{equation}
($x(t)\in\Rset^n$, $u(t)\in\Rset$, $y(t)\in\Rset$) and an objective function
\begin{equation}
z=\phi(y)
\end{equation}
to be optimized at steady state, i.e., we are looking for an output
feedback controller that steers the state to an operating point
$(x^{\star},u^{\star})$ that fulfills
\begin{equation}
\begin{array}{l}
 f(x^{\star},u^{\star})=0\\
\phi(h(x)) \mbox{ is locally max. w.r.t.\ $u\in\Rset$ at } x^{\star}
\end{array}
\end{equation}
where the functions $f$ and $\phi$ are unknown or partially
known.\\

We follow the statement of the general nonlinear problem of extremum seeking 
given in \cite{AK03}, that we adapt here to an output feedback framework.

{\bf Assumption A1 (Existence of stabilizing output feedback).} There
exists a smooth output feedback $u(t)=\alpha(y(t),\bar u, \bar y)$
parametrized by the parameter pair $(\bar u,\bar y)\in{\cal
  U}\times{\cal Y}$ (a pair of reference input and output), such that
\begin{equation}
\label{noninvasive}
\alpha(\bar y,\bar u,\bar y)=\bar u \ ,
\end{equation}
and the closed loop system
\begin{equation}
\label{closedloop}
\dot x = f(x,\alpha(h(x),\bar u,\bar y))
\end{equation}
admits a unique equilibrium $x_{eq}(\bar u, \bar y)$, that is locally
asymptotically stable for any $(\bar u, \bar y)\in {\cal U}\times{\cal
  Y}$.  

For example, for the chemostat
\eqref{chemostat},\,\eqref{eq:output} the output feedback could be of
the form $u(t)=\bar u+G_1(\bar s-s(t))$ with a sufficiently large
$G_1$ \cite{SRRD13}.

Under Assumption A1, we consider the extremum-seeking problem for the
closed-loop dynamics (\ref{closedloop}) as if the pair $(\bar u,\bar y)$ was a
new control, and look for pairs such that $h(x_{eq}(\bar u,\bar
  y))=\bar y$. For this purpose, we shall consider the output-input characteristics defined as follows.\\

{\bf Assumption A2 (Output-input characteristics).} There exists a smooth function $\psi: {\cal
  Y}\mapsto {\cal U}$ such that
\begin{equation}\label{eq:assA2}
\bar y = h(x_{eq}(\bar u, \bar y)) \; \Leftrightarrow \;
\bar u=\psi(\bar y)
\end{equation}
for any $\bar y \in {\cal Y}$.

This assumption means that for each parameter $\bar y \in {\cal
    Y}$, there exists an unique input $\bar u \in {\cal U}$ that is a zero of the
    function $u \mapsto h(x_{eq}(u,\bar y))-\bar y$. For instance, in
    the chemostat model, positive equilibriums have to fulfill 
    $\alpha(s_{eq},\bar u,\bar s)=\mu(s_{eq})$, and then having
    $s_{eq}=\bar s$, along with (\ref{noninvasive}), amounts to write $\bar u=\mu(\bar y)$ .\\

{\bf Assumption A3 (Existence of local maximum).} There exists $y^{\star} \in {\cal Y}$ such
that the function $\phi(\cdot)$ has a local strict maximum at
$y^{\star}$.\\

Finally, we shall assume that the unknown objective function
$\phi(\cdot)$ possesses a known structure.\\

{\bf Assumption A4 (Structure of objective function).} There exists a smooth function $\varphi: \Rset^{2} \mapsto
\Rset$
such that 
\begin{equation}
\phi(y)=\varphi(y,\psi(y)) \ , \forall y \in {\cal Y}
\end{equation}

Assumptions A2 and A4 imply that the extremum seeking problem amouts to
optimize a known function of $y$ and $u$ at steady state. Thus, the
problem consists then in finding an output
feedback strategy that drives the state of the system to $x^{\star}$
such that $h(x^{\star})=y^{\star}$ using the knowledge of the feedback
law $\alpha$ and the objective $\varphi$ only. The following two
sections present two approaches.

\section{Continuous-time adaptation}
\label{sec:cont}
The continuous-time approach first constructs a differential equation
for $\bar u$ that achieves $\bar u=\psi(\bar y)$ (and, thus, $\bar
y=h(\bar x)$) asymptotically, and then applies a gradient search
(again through a differential equation) along the curve of $(\bar
u,\bar y)$ given implicitly by $\bar y=h(\bar x(\bar u,\bar y))$. Note
that for feedback laws of the form $u(t)=\bar u+G[\bar y-y(t)]$ this
implicit curve corresponds to the curve of open-loop equilibrium
outputs.
\subsection{Step 1: an adaptive scheme for $\bar u$}
\label{sec:adaptbaru}
We look for a dynamics
\begin{equation}
\dot{\bar u} = \beta(y,\bar u,\bar y)
\end{equation}
with $\beta(\bar y,\bar u,\bar y)= 0$, such that
\begin{equation}
\bar E(\bar y)=\left[\begin{array}{c}
x_{eq}(\psi(\bar y),\bar y)\\
\psi(\bar y)
\end{array}\right]
\end{equation}
is a  locally asymptotically stable equilibrium of the coupled dynamics
\begin{equation}
\label{coupled-adaptive}
\begin{array}{lll}
\dot x & = & f(x,\alpha(h(x),(\bar u,\bar y)))\\
\dot{\bar u} & = & \beta(h(x),\bar u,\bar y)
\end{array}
\end{equation}
Note that Assumptions A1 and A2 ensures that $E(\bar y)$ is an
equilibrium of (\ref{coupled-adaptive}) for any $\bar y \in {\cal Y}$,
and that one has
\begin{equation}
\lim_{t\to +\infty} y(t)=\bar y \ .
\end{equation} 
This is a classical adaptive output control problem, for which several
techniques are available in the literature \cite{AW94,KKK95}.  We
simply require the convergence
\begin{equation}
\lim_{t\to +\infty} \bar u(t)=\psi(\bar y)
\end{equation} 
to be uniform w.r.t. $\bar y \in {\cal Y}$.

\subsection{Step 2: a continuation method for $\bar y$}

Step 1 with Assumption A4 provides an approximation of
the unknown objective function $\phi(\cdot)$ at a fixed $\bar y$:
\begin{equation}
\phi(\bar y) \simeq \bar \phi(\bar y)=\varphi(\bar y,\bar u) \ .
\end{equation}
The continuation consists in having $\bar y$ evolving slowly
(step-wise or continuously) with a steepest descent based on an
estimation of the gradient of the function $\phi$.
A continuous adaptation can be written as
\begin{equation}
\dot{\bar y}=-\epsilon\left[\partial_{1}\varphi(\bar y,\bar u)
+\partial_{2}\varphi(\bar y,\bar u)\hat\psi^{\prime}(\bar y,\bar u)\right]
\end{equation}
where $\epsilon>0$ is small compared to the time scale of dynamics
(\ref{coupled-adaptive}), and
$\hat\psi^{\prime}(\bar y,\bar u)$ is an estimation of the derivative
of the function $\psi$ at $\bar y$, for which $\bar u$ is an
estimation of $\psi(\bar y)$.
Several gradient estimations techniques are available in the
literature \cite{ME03,MCB09,YIKLL10,FSB12,BFB12}.
In the scalar case, an estimation of the sign of the derivative of the
function $\phi$ is enough to choose the steepest descent. 
We can use for instance a dynamics with delay:
\begin{equation}
\begin{array}{lll}
\frac{1}{\epsilon}\dot{\bar y} & = & \delta(\bar y,\bar u,\bar y(t-\tau),\bar u(t-\tau))\\
&=&  \sgn\left[(\varphi(\bar y,\bar u)-\varphi(\bar
  y(t-\tau),u(t-\tau))(\bar y-\bar y(t-\tau))\right] \ .
\end{array}
\end{equation}
Then, the overall dynamics
\begin{equation}
\begin{array}{rll}
\dot x & = & f(x,\alpha(h(x),\bar u,\bar y))\\
\dot{\bar u} & = & \beta(h(x),\bar u,\bar y)\\
\dot{\bar y} & = & \epsilon \delta(\bar y,\bar u,\bar
y(t-\tau),\bar u(t-\tau))
\end{array}
\end{equation}
is slow-fast where $\bar u=\psi(\bar y)$ is the (attracting) critical
manifold.\\

Finally, the output feedback strategy that we consider takes the
following form
\begin{equation}
u=\alpha(y,\bar u,\bar y) \mbox{ with }
\left\{\begin{array}{l}
\dot{\bar u} = \beta(y,\bar u,\bar y)\\
\dot{\bar y} = \epsilon \delta(\bar y,\bar u,\bar
y(t-\tau),\bar u(t-\tau))
\end{array}\right.
\end{equation}
Let us underline that the output $y$ is naturally filtered by the dynamics to obtain the
estimation $y^{\star}$ as
\begin{equation}
\hat{y}^{\star}=\lim_{t\to+\infty} \bar y(t)
\end{equation}
that provides a robustness w.r.t. to measurement noise. 
Of course, the price to pay is to have a slow convergence due to small $\epsilon$.

\section{Extremum-seeking using numerical optimization algorithms}
\label{sec:numerics}
We compare the continuous-time adaptation with classical numerical
optimization algorithms that are difficult to express as dynamical
systems such as golden-section iteration (which is applicable for
single-parameter problems) and Newton iterations. These algorithms do
not require any knowledge of the underlying system apart from the
ability to evaluate the asymptotic output $\lim_{t\to\infty}y(t)$ for
any given admissible inputs $(\bar u,\bar y)$. Once the stabilizing
feedback law $u(t)=\alpha(y(t),\bar u,\bar y)$ is implemented, the
problem is reduced to a classical numerical optimization problem:
$\phi(h(x_{eq}(\bar u, \bar y)))\to\max$. In fact, if the feedback law
is linear of the type
\begin{equation}
  \label{eq:feedbackgen}
  u(t)=\alpha(y(t),\bar u,\bar y)=\bar u+G[\bar y-y(t)]
\end{equation}
this is an optimization problem in the single parameter $\bar v=\bar
u+G\bar y$:
\begin{equation}
\phi(h(x_{eq}(\bar u+G\bar y)))=:F(\bar v)\to\max\mbox{\quad w.r.t
$\bar v$.}\label{eq:optimization}
\end{equation}
We test a combination of two classical optimization algorithms. We
start with an initial admissible interval $[\bar v_\mathrm{low},\bar
v_\mathrm{up}]$. Then we use golden section search to narrow the
initial interval down to a given fraction of its original size and
then apply a Newton iteration to solve $F'(\bar v)=0$ in the remaining
interval. Neither of these two methods requires knowledge about the
structure of $F$, only the ability to evaluate $F$ and its derivatives
at any desired point (where derivatives can be approximated by finite
differences). For the evaluation of $F$ and its derivatives in a given
point $\bar v_0$ we use a naive \emph{act-and-wait} approach
\cite{I06}: one sets the parameter $\bar v=\bar v_0$, such that the
feedback control law is $u(t)=\bar v_0-Gy(t)$, waits until the
transients have settled (such that $x=\bar x_0:=x_{eq}(\bar v_0)$),
and then reads off $y=h(\bar x_0)$ and evaluates 
$\phi(y)=\varphi(y,\bar v_{0}-Gy)$. For the
golden section search we iterate from an initial interval $[\bar
v_1,\bar v_3]$ by setting $\bar v_2=\bar v_1+(\bar v_3-\bar
v_1)(3-\sqrt{5})/2$ and evaluating $F_1=F(\bar v_1)$, $F_2=F(\bar
v_2)$ and $F_3=F(\bar v_3)$ to obtain an initial triplet $(\bar v_1,
\bar v_2, \bar v_3)$ and the corresponding objective function values
$(F_1,F_2,F_3)$. Then we proceed with the standard golden search to
iteratively obtain new triplets $(\bar v_1, \bar v_2, \bar v_3)$ until
$\bar v_3-\bar v_1<\mathtt{tol}$.  Correspondingly, for the Newton
iteration with an initial guess $\bar v_\mathrm{old}$ we evaluate $F$
at the three points $(\bar v_\mathrm{old}-h,\bar v_\mathrm{old},\bar
v_\mathrm{old}+h)$ (where $h$ is a small finite-difference deviation)
and choose as the new point $\bar v_\mathrm{new}$ the maximum point of
the interpolating parabola. The golden section search is known to
converge globally for unimodal functions with rate $(\sqrt{5}-1)/2$,
whereas the Newton iteration converges quadratically close to a local
maximum. In the presence of random disturbances in the evaluation of
$F$ we gradually increase the accuracy of the evaluation of $F$ during
the iteration by using the mean of $y$ over a longer time interval
after the transients have settled. Note that, when using act-and-wait
and a discrete-time optimization algorithm one has removed two of the
two slower time scales present in \cite{AK04}, replacing them with the
algorithmic iteration (the Newton iteration converges
super-exponentially).

\section{Extremum seeking for the chemostat}
\label{section_chemostat}
Reference~\cite{SRRD13} has shown that there exists an output feedback
that satisfies the assumptions A1, A2 and A4 for the model
(\ref{chemostat}), with a simple saturated proportionate controller.
Define first the saturation function as follows.
\begin{equation}
  \sat_{[D_{\min},D_{\max}]}(\xi)=\left|\begin{array}{ll}
    D_{\max} &\mbox{if $\xi> D_{\max}$,}\\
    \xi &\mbox{if $\xi\in[D_{\min}, D_{\max}]$,}\\
    D_{\min} &\mbox{if $\xi< D_{\min}$,}
  \end{array}\right.
\end{equation}
For this saturation function \cite{SRRD13} proved the following
proposition.
\begin{prop}\label{thm:feedback:stable}
  Assume that the reference values $(\bar s,\bar D)$ are in
  $[s_{\min},s_\mathrm{in})\times[D_{\min},D_{\max}]$, where the bounds satisfy
\begin{equation}
\begin{array}{ll}
  D_{\min}&< \mu(s)\mbox{\quad for all $s\in[s_{\min},s_\mathrm{in}]$,}\\
  D_{\max}&>\mu(s)\mbox{\quad for all $s\in[0,s_\mathrm{in}]$.}
\end{array}
\end{equation}
Then the output feedback
\begin{equation}\label{chemostat-feedback}
  u(y,\bar D,\bar s)=\sat_{[D_{\min},D_{\max}]} \left(\bar D
  - G_{1}(y-\bar s)\right)
\end{equation}
with a gain $G_{1}$ satisfying
\begin{equation}
G_{1} > \max \left( \max_{s\in[0,s_\mathrm{in}]}-\mu'(s) \; , \;
\frac{\bar D-\mu(s_\mathrm{in})}{s_\mathrm{in}-\bar s} \right)
\end{equation}
guarantees that the closed-loop dynamics \eqref{chemostat},
\eqref{chemostat-feedback} has a stable equilibrium
$(s_\mathrm{eq},b_\mathrm{eq})\in[0,s_\mathrm{in})\times(0,\infty)$,
which attracts all initial conditions
$(s(0),b(0))\in[0,s_\mathrm{in})\times(0,\infty)$.
\end{prop}
See \cite{SRRD13} for a proof.\\

From system \eqref{chemostat}, one has
at steady state
\begin{equation}
\mu(s_\mathrm{eq})=\bar D -G_{1}(s_\mathrm{eq}-\bar s)\mbox{,}
\end{equation}
and consequently one has
\begin{equation}
s_\mathrm{eq}=\bar s \; \Leftrightarrow \; \bar D =\mu(s_\mathrm{eq})\mbox{,}
\end{equation}
that is, Assumption A2 is fulfilled.\\

Note that in system \eqref{chemostat} at steady state one has
$b=s_\mathrm{in}-s$. Thus, the objective function
\begin{equation}
\phi(s)=\mu(s)(s_\mathrm{in}-s)\mbox{,}
\end{equation}
which is non-negative such that $\phi(0)=\phi(s_{in})=0$, is
equivalent to the production rate $\mu(s)b$. For $\phi$ the
assumptions A3 and A4 are fulfilled.

\subsection{Implementation of continuous-time adaptation}
\label{sec:cont:chemostat}
For the continuous-time adaptation we can choose the same adaptive
scheme for $\bar u$ for the step~1 described in
Section~\ref{sec:adaptbaru} as was chosen in \cite{SRRD13} for
the identification of $\mu$.
\begin{prop}\label{thm:bard:adapt}
For any $\bar s \in (0,s_{in})$ and numbers $D_{\min}$, $D_{\max}$ 
such that $0<D_{\min} < \mu(\bar s) < D_{\max}$, the dynamical feedback law
\begin{equation}
\label{dyn-feedback}
\begin{array}{rll}
u(y,\bar D,\bar s) & =  & \sat_{[D_{\min},D_{\max}]}\left(\bar D - G_{1}(y-\bar s)\right)\\
\dot{\bar D} & = & -G_{2}(y-\bar s)(\bar D-D_{\min})(D_{\max}-\bar D)
\end{array}
\end{equation}
exponentially stabilizes the system \eqref{chemostat} locally about
$(\bar s,s_{in}-\bar s)$, for any positive constants
$(G_{1},G_{2})$ such that $G_{1}> -\mu^{\prime}(\bar s)$.  Furthermore
one has 
\[
\lim_{t\to+\infty}\bar D(t)=\mu(\bar s)
\]
\end{prop}
See \cite{SRRD13} for the proof.\\

Our continuous-time extremum seeking output feedback for the chemostat
model \eqref{chemostat} is thus given by the following equations
\begin{align}
\label{ext_seek_chemostat_u}
  u(y,\bar D,\bar s)  = &\ \sat_{[D_{\min},D_{\max}]}\left(\bar D - G_{1}(y-\bar s)\right)\\
\label{ext_seek_chemostat_bard}
  \dot{\bar D} = &\ -G_{2}(y-\bar s)(\bar D-D_{\min})(D_{\max}-\bar
  D)\\
\label{ext_seek_chemostat_bars}
  \frac{1}{\epsilon}\dot{\bar s} =&\ \sgn\left[
    (\bar D(s_{in}-\bar s)-\bar D(t-\tau)(s_{in}-\bar s(t-\tau)))(\bar s-\bar
    s(t-\tau))\right]
\end{align}

\subsection{Implementation of discrete-time optimization}
\label{sec:opt:chemostat}
The discrete-time optimization uses only the saturated feedback, combining the parameters $\bar D$ and
$\bar S$ into a single parameter $\bar v=\bar D+G_1\bar s$, such that
\begin{equation}
  \label{eq:feedback:single}
  u(y,\bar v)  =  \sat_{[D_{\min},D_{\max}]}(\bar v - G_{1}y)\mbox{.}
\end{equation}
If the initial bracketing interval $[\bar v_1,\bar v_2]$ for $\bar v$
is chosen such the line $\ell(\bar v)=\{(\bar D,\bar s): \bar
D+G_1\bar s=\bar v\}$ intersects the domain of admissible reference
values $[s_{\min},s_\mathrm{in})\times[D_{\min},D_{\max}]$ for all
$\bar v\in[\bar v_1,\bar v_2]$, Proposition~\ref{thm:feedback:stable}
ensures the existence of a unique equilibrium $s_\mathrm{eq}(\bar v)$
for all $\bar v\in[\bar v_1,\bar v_2]$. Thus, we can proceed with the
golden section search for the objective functional
\begin{equation}
  \label{eq:chem:obj}
  F(\bar v)=u(s_\mathrm{eq}(\bar v),\bar v)[s_\mathrm{in}-s_\mathrm{eq}(\bar v)]\mbox{,}
\end{equation}
where $u$, as defined in \eqref{eq:feedback:single}, is the asymptotic
value of the input $u$.  Once, the golden section search has shrunk
the bracketing interval to a certain size, we apply a single step of
the Newton iteration for the objective functional $F$.
\section{Numerical simulations}
\label{section_demo}
The output feedback law
\eqref{ext_seek_chemostat_u}--\eqref{ext_seek_chemostat_bars} has been
simulated for a Haldane function
\begin{equation}
\label{Haldane}
\mu(s)=\frac{\mu_{\max}s}{K+s+s^{2}/K_{i}}
\end{equation}
and the following values of the parameters:
\begin{center}
\begin{tabular}{|cccccc|}
\hline
$\mu_{\max}$ & $K$ & $K_{i}$ & $S_{in}$ & $D_{\min}$ & $D_{\max}$\\
\hline
1.0 & 1.0 & 0.1 & 1.0 & 0.0 & 1.0\\
\hline
\end{tabular}
\end{center}
Note however, that in our simulations we make no assumptions about the
type of $\mu$ apart from those required in the assumptions of
propositions \ref{thm:feedback:stable} and \ref{thm:bard:adapt}
(leading to the validity of assumptions A1--A4). For the control laws
in \eqref{ext_seek_chemostat_u}--\eqref{ext_seek_chemostat_bars}, the
following parameters have been chosen:
\begin{center}
\begin{tabular}{|ccc|}
\hline
$G_{1}$ & $G_{2}$ & $\epsilon$\\
\hline
1.0 & 1.0 & $1e^{-3}$\\
\hline
\end{tabular}
\end{center}
To simulate measurement noise, we have considered that the output is
corrupted in the following way
\begin{equation}
y(t)=s(t)(1+b(t))\label{eq:disturbance}
\end{equation}
where $b(\cdot)$ is a random square signal of frequency $\omega$ and
amplitude $a$, whose values are given below.
\begin{center}
\begin{tabular}{|cc|}
\hline
$\omega$ & $a$\\
\hline
0.2 & 0.05\\
\hline
\end{tabular}
\end{center}

\subsection{Continuous-time adaptation}
\label{sec:cont:results}
\begin{figure}[ht]
  \centering
  \includegraphics[width=7cm]{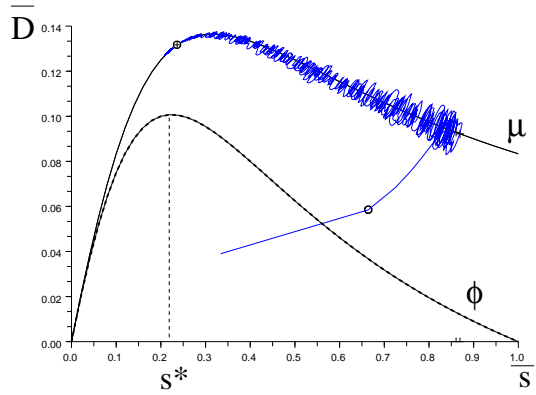}
  \caption{Trajectory of \eqref{chemostat} with feedback laws
    \eqref{ext_seek_chemostat_u}--\eqref{ext_seek_chemostat_bars} in the
    $(\bar s,\bar D)$ plane.}
  \label{fig:Ds}
\end{figure}
\begin{figure}[ht]
  \centering
  \includegraphics[width=7cm]{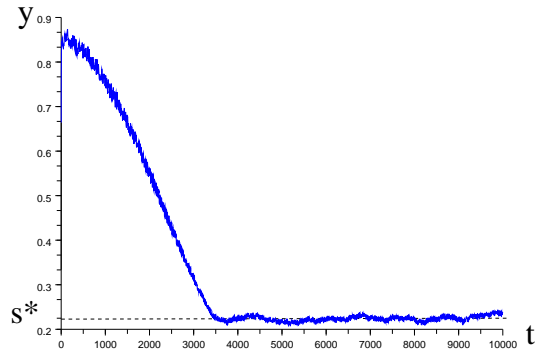}
  \caption{Output of \eqref{chemostat} with feedback laws
    \eqref{ext_seek_chemostat_u}--\eqref{ext_seek_chemostat_bars},
    disturbed by \eqref{eq:disturbance}, over time (same run as Figure~\ref{fig:Ds}).}
  \label{fig:ts}
\end{figure}
\begin{figure}[ht]
  \centering
  \includegraphics[width=7cm]{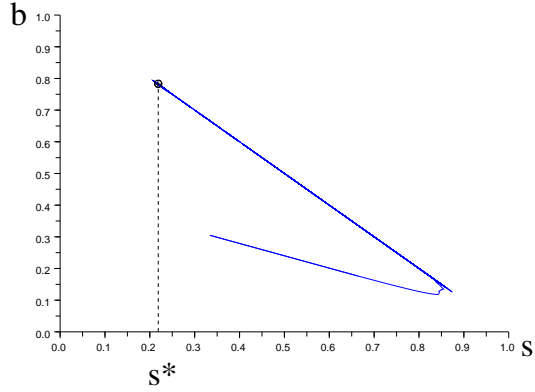}
  \caption{Trajectory of same run as Figures~\ref{fig:Ds} and
    \ref{fig:ts} in the $(s,b)$ plane ($\{(s,b): s+b=s_\mathrm{in}\}$
    is invariant).}
  \label{fig:sb}
\end{figure}
Figure~\ref{fig:Ds} shows how the trajectory of the variables $(\bar
s,\bar D)$ converges to the graph of the function $\mu(\cdot)$ i nthe
$(\bar s,\bar D)$-plane initially up to disturbances due to the
measurement noise. This convergence toward $\mu(\cdot)$ is effected by
the control laws \eqref{ext_seek_chemostat_u} and
\eqref{ext_seek_chemostat_bard} on the fast timescale (for $t\ll 100$,
not visible in the time series depicted in Figure~\ref{fig:ts}).  Once
the feedback laws \eqref{ext_seek_chemostat_u} and
\eqref{ext_seek_chemostat_bard} have achieved quasi-equilibrium for
fixed $\bar s$ the dynamics of \eqref{ext_seek_chemostat_bars}
maximizes the functional $\phi$ on the slow timescale. Timescale
separation $\epsilon=10^{-3}$ is chosen for our simulations. In
Figures~\ref{fig:Ds} and \ref{fig:sb} one can also see that the
adaptation of $\bar s$ via \eqref{ext_seek_chemostat_bars} forces the
state $(s,b)$ to converge near the equilibrium $(s^{\star},b^{\star})$
(indicated as intersection between dotted line and graph $\mu(\cdot)$,
where $s^{\star}$ is maximizing the function $\phi(\cdot)$. We note
that the adaptation has converged already after $t=3500$.

\emph{Remark} If we choose the timescales not sufficiently different
(say, $\epsilon=10^{-2}$ in our example) we observe that the chemostat
\eqref{chemostat} with feedback laws
\eqref{ext_seek_chemostat_u}--\eqref{ext_seek_chemostat_bars} becomes
dynamically unstable in in its fixed point $(s,b,\bar s,\bar
D)=(s^{\star},b^{\star},s^{\star},\mu(s^\star))$.
\subsection{Discrete-time algorithms}
\label{sec:disc:results}
\begin{figure}[ht]
  \centering
  \includegraphics[width=10cm]{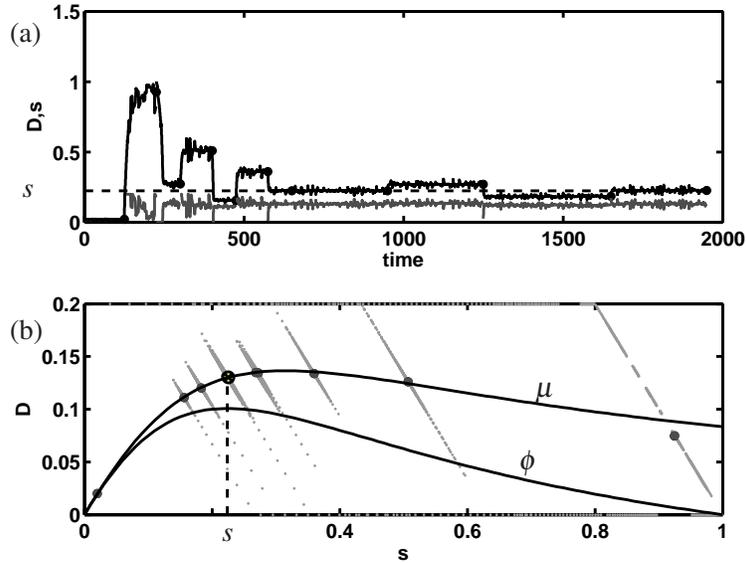}
  \caption{Time profiles (a) and input-output plane (b) for
    time-discrete optimization. The grey time profile in (a) is
    $u(t)=D(t)$ and the black profile is $s(t)$. Full circles in (a) and
    (b) indicate inputs and outputs where the objective $F$ was
    evaluated.}
  \label{fig:disc}
\end{figure}
For the discrete-time optimization algorithm we chose an initial
bracketing interval $[\bar v_1, \bar v_2]=[0.04,1]$ (gain $G_1=1$ in
feedback law \eqref{eq:feedback:single} identical to the
continuous-time adaptation). During the golden section search (until
$t\approx700$ in Figure~\ref{fig:disc}(a)) we checked every
$t_\mathrm{inc}=25$ time units if the transients have settled
(criterion is a decrease of the standard deviation compared to the
previous time interval). Once the output is accepted as stationary we
used the mean of output and input over the last time interval to
calculate the objective $F$. A larger $t_\mathrm{inc}$ inceases the
accuracy of the evaluation of $F$ in the presence of disturbances due
to averaging. We observe that after $7$ evaluations of $F$ the golden
section search reduces the bracketing interval to length $0.2$ and the
estimated optimal value of $s$ is already close to $s^\star$
($t\approx700$ in Figure~\ref{fig:disc}(a)).  The Newton step then
only confirms the optimality of the final estimate of $s^\star$ (up to
disturbance level). For the Newton iteration we chose $t_\mathrm{inc}$
larger ($t_\mathrm{inc}=100$) to increase the accuracy of $F$
evaluations and use the points $(\bar v-h,\bar v,\bar v+h)$ with
$h=0.05$ to evaluate the approximating parabola for $F$ in $\bar
v$. Figure~\ref{fig:disc}(b) shows the trace (grey dots) of outputs
$s$ and inputs $u$ (equaling dilution rate $D$) in the
$(s,D)$-plane. The gaps in the trace are due to finite sampling
time. The full circles show the points where the transients were
accepted as settled.

The speed and accuracy of the discrete-time optimization is limited by
the choice of $t_\mathrm{inc}$, which in turn has to be determined
sufficiently large to enable an averaging effect for the disturbances.
\section{Conclusion}
We have proposed two extremum-seeking schemes: one approach uses
dynamical output feedback that is based on a continuation method.  Its
closed-loop dynamics possesses two time scales: the timescale of the
original system with stabilizing feedback and a slow timescale for the
adaptation scheme. As an alternative we test classical optimization
algorithms, applied to the steady-state outputs of the
feedback-stabilized system, which have a continuous timescale (which
is exponentially convergent) and a discrete timescale (which is
exponentially or superexponentially convergent). This is in contrast
to extremum seeking techniques that use dither signals, which need
three timescales for nonlinear systems (the time scale of the
stabilized system is the fastest). We illustrate our approaches on the
chemostat model for unknown growth function and single measurement.
The numerical simulations suggest that the method is practically
viable and robust with respect to measurement disturbances. Extensions
and a detailed convergence analysis for general systems and the
continuous-time adaptation will be the subject of future
work. Conditions for the convergence of the discrete-time optimization
can be deduced directly from the corresponding statements for the
original numerical algorithms.

\end{document}